\documentclass[11pt,makeidx]{article}
\input amssym.def
\input amssym.tex
\usepackage{epsfig,psfig}
\newtheorem{theorem}{Theorem}
\newtheorem{lemma}{Lemma}

\newtheorem{coro}{Corollary}

\def\biblitem#1{\bibitem{#1}}

\def\slfrac#1#2{\hbox{\kern.1em %
 \raise.5ex\hbox{\the\scriptfont0 #1}\kern-.11em %
 /\kern-.15em\lower.25ex\hbox{\the\scriptfont0 #2}}}

\def\pmod#1{\allowbreak \mkern9mu ({\rm mod}\,\,#1)}
\def\legendre#1#2{\left(#1\over #2\right)}
\newcommand{\eqn}[1]{(\ref{#1})}
\newcommand{\hsp}{\hspace*{\parindent}}

\newcommand{\eeq}{\end{equation}}
\newcommand{\beql}[1]{\begin{equation}\label{#1}}
\newcommand{\bsq}{{\vrule height .9ex width .8ex depth -.1ex }}

\newcommand{\CN}{C^{(N)}}
\newcommand{\LAP}{\Lambda^{\ast_{\Pi}}}
\newcommand{\apo}{{\ast_{\Pi_1}}}
\newcommand{\apt}{{\ast_{\Pi_2}}}
\newcommand{\LA}{\Lambda^\ast}

\newcommand{\sM}{{\cal M}}
\newcommand{\ZZ}{{\Bbb Z}}
\newcommand{\RR}{{\Bbb R}}
\newcommand{\FF}{{\Bbb F}}
\newcommand{\QQ}{{\Bbb Q}}

\newcommand{\om}{\omega}

\newcommand{\sg}{\sigma}
\newcommand{\La}{\Lambda}
\newcommand{\Ga}{\Gamma}

\newcommand{\Om}{\Omega}

\newcommand{\sL}{{\cal L}}

\makeatletter
\def\@sect#1#2#3#4#5#6[#7]#8{\ifnum #2>\c@secnumdepth
     \def\@svsec{}\else
     \refstepcounter{#1}\edef\@svsec{\csname the#1\endcsname.\hskip .75em }\fi
     \@tempskipa #5\relax
      \ifdim \@tempskipa>\z@
        \begingroup #6\relax
          \@hangfrom{\hskip #3\relax\@svsec}{\interlinepenalty \@M #8\par}%
        \endgroup
       \csname #1mark\endcsname{#7}\addcontentsline
         {toc}{#1}{\ifnum #2>\c@secnumdepth \else
                      \protect\numberline{\csname the#1\endcsname}\fi
                    #7}\else
        \def\@svsechd{#6\hskip #3\@svsec #8\csname #1mark\endcsname
                      {#7}\addcontentsline
                           {toc}{#1}{\ifnum #2>\c@secnumdepth \else
                             \protect\numberline{\csname the#1\endcsname}\fi
                       #7}}\fi
     \@xsect{#5}}
\def\@begintheorem#1#2{\it \trivlist \item[\hskip \labelsep{\bf #1\ #2.}]}

\def\plain{plain}\ifx\fmtname\plain\csname fi\endcsname
     
     \input docstrip
     \preamble

     Do not distribute the stripped version of this file.
     The checksum in the header refers to the documented version.

     \endpreamble
     \generateFile{here.sty}{t}{\from{here.doc}{}}
     \endinput
\fi
\ifcat a\noexpand @\let\next\relax\else\def\next{%
    \documentstyle[here,doc]{article}\MakePercentIgnore}\fi\next
\ifx\@Hxfloat\@Hundef\else\expandafter\endinput\fi
\let\@Hxfloat\@xfloat
\def\@xfloat#1[{\@ifnextchar{H}{\@HHfloat{#1}[}{\@Hxfloat{#1}[}}
\def\@HHfloat#1[H]{%
\expandafter\let\csname end#1\endcsname\end@Hfloat
\vskip\intextsep\vbox\bgroup\def\@captype{#1}\parindent\z@
\ignorespaces}
\def\end@Hfloat{\egroup\vskip \intextsep}

\makeatother

\begin{document}
\begin{center}
{\Large {\bf The Shadow Theory of Modular and Unimodular Lattices}} \\
\vspace{1.5\baselineskip}
{\em E. M. Rains} and {\em N. J. A. Sloane} \\
\vspace*{1\baselineskip}
Information Sciences Research, AT\&T Labs-Research \\
180 Park Avenue, Florham Park, NJ 07932-0971, U.S.A. \\
\vspace{1.5\baselineskip}
Apr. 27, 1998 \\
\vspace{1.5\baselineskip}
{\bf ABSTRACT}
\vspace{.5\baselineskip}
\end{center}
\setlength{\baselineskip}{1.5\baselineskip}

It is shown that an $n$-dimensional unimodular lattice has minimal norm at most
$2[n/24] +2$, unless $n=23$ when the bound must be increased by 1.
This result was previously known only for even unimodular lattices.
Quebbemann had extended the bound for even unimodular lattices to strongly
$N$-modular even lattices for $N$ in
$$
\{1,2,3,5,6,7,11,14,15,23 \} ~,
\eqno{(\ast)}
$$
and analogous bounds are established here for odd lattices satisfying
certain technical conditions (which are trivial for $N = 1$ and $2$).
For $N> 1$ in $(\ast)$, lattices meeting the new bound are constructed that are analogous to the
``shorter'' and ``odd'' Leech lattices.
These include an odd associate of the 16-dimensional Barnes-Wall
lattice and shorter and odd associates of the Coxeter-Todd lattice.
A uniform construction is given for the (even) analogues of the
Leech lattice,
inspired by the fact that $(\ast)$ is also the set of square-free orders
of elements of the Mathieu group $M_{23}$.
\clearpage
\thispagestyle{empty}
\setcounter{page}{1}

\section{Introduction}
\hsp
The study of unimodular lattices (i.e. integral lattices of determinant 1) is an important chapter in the classical theory of quadratic forms.
Another way to characterize a unimodular lattice is
that it is equal to its dual.
A {\em modular} lattice (the term was introduced by Quebbemann \cite{Q95};
see also \cite{Q97},
\cite{Q98}) is an integral lattice which is geometrically similar to its dual.

In other words, an $n$-dimensional integral lattice $\La$ is {\em modular}
if there exists a similarity $\sg$ of $\RR^n$ such that
$\sg (\La^\ast ) = \La$, where $\La^\ast$ is the dual lattice.
If $\sg$ multiplies norms by $N$, $\La$ is said to be
{\em $N$-modular}.
For example, the sporadic root lattices $E_8$,
$F_4 ( \cong D_4 )$, $G_2 ( \cong A_2 )$ are respectively 1-, 2- and 3-modular.
In the last two cases the modularity maps short roots to long roots.

If $N$ is a composite number, a {\em strongly $N$-modular} lattice \cite{Q97} satisfies certain additional conditions given in Section 3.

To date the study of $N$-modular lattices for $N > 1$ has focused
on even lattices, but in the present paper we remove this restriction
and also consider odd lattices.

The simplest example of an $N$-modular lattice for $N$ prime is the two-dimensional lattice
$C^{(N)} = \ZZ \oplus \sqrt{N} \ZZ$.
The similarity $\sg$ takes $(x,y)$ to
$(\sqrt{N} y, \sqrt{N} x )$, and maps $C^{(N) \ast}$ to
$C^{(N)}$.
More generally, for any positive integer $N$,
$$\CN = \sum_{d|N} ~ \sqrt{d} \ZZ$$
is a strongly $N$-modular lattice of dimension equal to $d(N)$, the number of
divisors of $N$.

The main goal of this paper is to prove Theorems~\ref{th1} and \ref{th2}.
\begin{theorem}
\label{th1}
An $n$-dimensional unimodular lattice has minimal norm
\beql{Eq1}
\mu \le 2 \left[\frac{n}{24} \right] +2 ~,
\eeq
unless $n=23$ when $\mu \le 3$.
\end{theorem}

\paragraph{Remarks.}
(1)~The form of \eqn{Eq1} suggests that dimension 24 may be special, and of
course it is: there is a unique 24-dimensional lattice meeting the
bound, the Leech lattice $\La_{24}$ (cf. \cite{SPLAG}).
The best odd lattice in dimension 24 is the ``odd Leech lattice''
$O_{24}$ of minimal norm 3, and the
exception to the bound in dimension 23 is necessary because of the existence of the ``shorter Leech lattice''
$O_{23}$,
which also has minimal norm 3.

(2)~Theorem~\ref{th1} is the strongest upper bound presently known for unimodular lattices.
For {\em even} unimodular lattices this was already known \cite{Me41}, but for odd unimodular lattices it was known only that
$$\mu \le \left[ \frac{n+6}{10} \right]$$
for all sufficiently large $n$ \cite{Me157}.

(3)~For self-dual codes the situation is similar.
For doubly-even self-dual codes it was shown in
\cite{Me41},
\cite{Me30}
that the minimal distance $d$ of a code of length $n$ satisfies
\beql{Eq3}
d \le 4 \left[\frac{n}{24} \right] +4 ~,
\eeq
and for singly-even self-dual codes
$$d \le 2 \left[ \frac{n+6}{10} \right] ~,
$$
unless $n=2$, 8, 12, 22, 24, 32, 48 and 72 when the bound must
be increased by 2 \cite{Me158}.
The analogue of Theorem~\ref{th1} is given in \cite{RainsSB},
where it is shown that \eqn{Eq3} holds for all self-dual
codes, unless $n \equiv 22$ $(\bmod~24)$ when the upper bound
must be increased by 2.

So in the coding analogue to Theorem~\ref{th1} there are infinitely many exceptions, not just one.
However, it seems very likely that equality can hold in \eqn{Eq1}
and \eqn{Eq3}, and in the bounds of Theorem~\ref{th2}, for only finitely
many values of $n$ (compare \cite{Me41}).

(4)~In the coding analogue of Theorem~\ref{th1}, it can be shown that any self-dual code of length $n \equiv 0$ $(\bmod~24)$ meeting the bound in \eqn{Eq3}
must be doubly-even.
We conjecture that if $n \equiv 0$ $(\bmod~24)$ any unimodular lattice meeting
the bound of Theorem~\ref{th1} must be even, although we have so far
not succeeded in proving this.

(5)~Krasikov and Litsyn \cite{KrLit91} have recently shown that for
doubly-even self-dual codes of length $n$, where $n$ is large,
\eqn{Eq3} can be improved to
$$d \le 0.166315 \ldots n + o(n) , \quad n \to \infty ~.$$
No analogous result is known for even unimodular lattices.

(6)~Theorem~\ref{th1} is included in Theorem~\ref{th2}, but is stated separately because of the importance of the unimodular case.

For strongly $N$-modular lattices we will restrict our attention to values of $N$ from the set
\beql{Eq4}
\{1,2,3,5,6,7,11,14,15,23 \} ~,
\eeq
for which the corresponding {\em critical dimensions}
$D_N = 24 d(N)/ \prod_{p|N} (p+1)$ are respectively
\beql{Eq5}
\{24,16,12,8,8,6,4,4,4,2 \} ~.
\eeq
\begin{theorem}\label{th2}
For $N$ in \eqn{Eq4}, an $n$-dimensional strongly $N$-modular lattice which is
rationally equivalent to the direct sum of $n/ \dim \CN$ copies of $\CN$ has minimal norm
\beql{Eq6}
\mu \le 2 \left[ \frac{n}{D_N} \right] + 2 ~,
\eeq
unless $N$ is odd and $n=D_N - \dim \CN$ when
\beql{Eq7}
\mu \le 3 ~.
\eeq
\end{theorem}

\paragraph{Remarks.}
(1)~The form of \eqn{Eq6} suggests that dimension $D_N$ may be special, and indeed in each case there is a unique lattice in that dimension meeting
the bound (see Section 2).

(2)~We will say that an $n$-dimensional strongly $N$-modular lattice $\La$
that meets the appropriate bound from Theorems~\ref{th1} or \ref{th2} is {\em extremal}.
This definition agrees with the historical usage for even lattices,
but for odd unimodular lattices extremal has generally meant minimal norm $[n/8]+1$.
There are just 11 such lattices with the latter property (SPLAG, Chap.~19).
In view of Theorem~\ref{th1} the more uniform definition proposed here seems
preferable.
A lattice satisfying the hypothesis of Theorem~\ref{th2}
is {\em optimal} if it has the highest minimal norm of any such
lattice with the same $n$ and $N$.
An extremal lattice is {\em a priori} optimal.

(3)~We conjecture that any extremal lattice of dimension a multiple of $D_N$
must be even (compare Remark (4) above).

(4)~The bound of Theorem~\ref{th2}
for $N \ge 11$ is quite weak, even for moderate values of $n$.
If $N=23$, for example, extremal lattices almost
certainly do not exist in dimensions above 4.
(Of course the analogous bounds for even lattices
\cite{Q97} are also weak.)

Section 2 gives a number of examples, some of which (the odd versions of the
Barnes-Wall and Coxeter-Todd lattices, and the shorter Coxeter-Todd lattice, for instance) appear to be new.

In Section 3 we study certain Gauss sums $\gamma_\Pi (\La )$ associated with a lattice $\La$,
show how Atkin-Lehner involutions act on theta series, and define the concept
of strong modularity.
Section 4 studies the shadow of a lattice.
For example, Theorem~\ref{ThOdd} shows that the norm of every vector in the shadow of an odd lattice is congruent to (oddity $\La$)/4 modulo $2 \ZZ_2$.
In Section 5 it is shown that the theta series of a lattice and its shadow are (essentially) invariant under the action of a certain modular group
$\frac{1}{2} \Gamma_0 (4N)^+$.
The main result of this section is
Corollary~\ref{C3}.

Section 6 contains the proofs of Theorems~\ref{th1} and \ref{th2}
(which make use of Corollary~\ref{Co3} from Section 3, Eq.~\eqn{Eq14a}
from Section 4, and Theorem~\ref{C2} and Corollary~\ref{C3} from
Section 5), as well as some identities for modular functions that may
be of independent interest.

In Section 7 we briefly discuss bounds for $N$-modular
lattices not covered by Theorem~\ref{th2}.
In the Appendix we prove a general result about the nonexistence of
modular lattices in certain genera.  Among other things this
implies that any 7- or 23-modular lattice must
satisfy the hypothesis of Theorem~\ref{th2}.

\section{Examples of extremal modular lattices}
\hsp
Many examples of modular lattices meeting the bounds of Theorems~\ref{th1} and
\ref{th2} (and of the analogous bounds in Section 7) can be found for instance in
\cite{Bach97}, \cite{SPLAG}, \cite{Mar96}, \cite{Nebe24}, \cite{Nebe31},
\cite{NP}, \cite{Q87}, \cite{Q95}, \cite{Q97}.
Other examples will be constructed here.
Some nonexistence results are given in \cite{NV} and \cite{Schar94} (see also
\cite{SchV}, \cite{SchaV96}).

For unimodular lattices, the highest possible minimal norm is known for
dimensions $n \le 33$ and 40--48 \cite{Me157}, \cite{CS98}, and in this range the bound of
Theorem~\ref{th1} is achieved precisely for $n=8,12,14 -24,32$ and 40--48.

For $N=2$, lattices achieving the bound of Theorem~\ref{th2} are known
(see e.g. \cite{Bach97} and \cite{Q87}) in dimensions $n=4$, 8--16, 20, 24, 28, 32, 36, 40, 44, 48, and do not exist for $n=2$, 6, 18, 34; the
existence for $n=22$, 26, 30, 38, 42, 46 is open.

For $N=3$, lattices meeting the bound of Theorem~\ref{th2} are known
(cf. \cite{Bach97}, \cite{Nebe31} and the present paper) for $n=4$--12,
16--24, 28, 32, and do not exist for $n=2$, 14, 26 and 50.

Less is known for larger values of $N$, for which we refer the reader to the table in \cite{SN}.
(This table also has further information about many of the above lattices.)

We begin our discussion of specific constructions by noting the following generalization of a construction given in
\cite[Chap. 7, Theorem 26]{SPLAG} and \cite{Bach97}:
if $C$ is an additive (but not necessarily linear) trace self-dual\footnote{That is, self-dual with respect to the inner product $Tr(u \cdot \overline{v} )$ \cite{Me223}, \cite{Self}.}
code over $\FF_4$ of length $n$ and minimal distance $d$, then
``Construction A''\footnote{In other words, take the real form of the complex
lattice $\{ u \in \ZZ [\omega ]^n: u \bmod~2 \in C \}$, where $\omega$ is a primitive cube root of unity.}
produces a 3-modular lattice in dimension $2n$ with minimal norm
$\mu = \min \{4, d\}$.
If $C$ is even so is the lattice (and if $C$ is odd the shadow
of the lattice is obtained by lifting the shadow of the code).

Since all lattices arising in this way share the common sublattice
$(\sqrt{2} A_2 )^n$, they are rationally equivalent to $(C^{(3)} )^N$, where
$C^{(3)} = \ZZ \oplus \sqrt{3} \ZZ$ arises from the code $C$ with generator matrix $[1]$.
Thus these lattices all satisfy the hypothesis of Theorem~\ref{th2}.
In particular, the hexacode
(with $n=6$, $d=4$) \cite[p.~82]{SPLAG} gives rise to the Coxeter-Todd lattice
$K_{12}$.
There are two related additive self-dual codes, the shorter $(n=5, d=3)$ and
odd $(n=6, d=3)$ hexacodes \cite{Me223}, \cite{HoehnSD}, \cite{Self}.
The latter can be taken to be the additive code generated by all cyclic
shifts of $1 \om 1 000$.
Under Construction A these codes become the shorter and odd Coxeter-Todd lattices
$S^{(3)}$ and $O^{(3)}$ (see Theorem~\ref{ThCon1}).
Other examples of good additive codes over $\FF_4$ from \cite{Me223},
\cite{Self} lead to optimal 3-modular lattices in dimensions $n \le 22$, including possibly new lattices
in dimensions 14, 18 and 22.
Construction A applied to the dodecacode ($n=12$, $d=6$, \cite{Me223},
\cite{HoehnSD}, \cite{Self}) gives rise to a neighbor of Nebe's
24-dimensional extremal 3-modular lattice \cite{Nebe24}, which has minimal
norm 6 rather than 4.

As remarked above, for each $N$ in \eqn{Eq4} there is an especially interesting
extremal strongly $N$-modular even lattice $E^{(N)}$ in the critical
dimension $D_N$, having minimal norm 4.
There is also a $D_N$-dimensional strongly $N$-modular {\em odd lattice}
$O^{(N)}$
of minimal norm 3, and,
when $N=1$, 3, 5, 7, 11, a {\em shorter lattice} $S^{(N)}$
of dimension $D_N-1$ (if $N=1$) or $D_N -2$ (if $N > 1$),
also with minimal norm 3 (see Table~\ref{TaE}).
The even lattices are well known, see \cite{Q95}, \cite{Q97}.
It turns out that there is a uniform construction for all the above lattices
(except for $O^{(N)}$ when $N$ is even).
\begin{table}[htb]
\caption{Extremal strongly $N$-modular even lattice $E^{(N)}$ in the critical
dimension $D_N$ and its odd $(O^{(N)} )$ and shorter
$(S^{(N)} )$ associates in dimensions $D_N$ and $D_N$-$\dim \, \CN$ respectively.}
\label{TaE}
$$
\begin{array}{|c|c|c|c|c|} \hline
N & D_N & E^{(N)} & O^{(N)} & S^{(N)} \\ \hline
1 & 24 & \La_{24} & O^{(1)} = O_{24} & S^{(1)} = O_{23} \\ 
2 & 16 & BW_{16} & O^{(2)} & \mbox{---} \\ 
3 & 12 & K_{12} & O^{(3)} & S^{(3)} \\ 
5 & 8 & Q_8 (1) & O^{(5)} & Q_6 (4)^{+2} \\ 
6 & 8 & G_2 \otimes F_4 & O^{(6)} & \mbox{---} \\ 
7 & 6 & A_6^{(2)} & O^{(7)} & \left[ \begin{array}{rrrr}
3&1&0&-1 \\ 1&3&1&0 \\ 0&1&3&1 \\ -1&0&1&3
\end{array} \right] \\ 
11 & 4 & \left[ \begin{array}{rrrr}4&1&0&-2 \\
1&4&2&0 \\ 0&2&4&1 \\ -2&0&1&4
\end{array} \right] & \left[ \begin{array}{cc}
3 & 1 \\ 1 & 4 \end{array}\right] \oplus \left[ \begin{array}{cc}
3 & 1 \\ 1 & 4 \end{array}\right] & \left[ \begin{array}{cc}
3 & 1 \\ 1 & 4 \end{array}\right] \\ 
14 & 4 & \left[ \begin{array}{rrrr}4&1&0&-1 \\
1&4&1&0 \\ 0&1&4&1 \\ -1&0&1&4
\end{array} \right] & \left[ \begin{array}{cc}
3 & 1 \\ 1 & 5 \end{array}\right] \oplus \left[ \begin{array}{cc}
3 & 1 \\ 1 & 5 \end{array}\right] & \mbox{---} \\ 
15 & 4 & \left[ \begin{array}{rrrr}4&2&2&1 \\
2&4&1&2 \\ 2&1&6&3 \\ 1&2&3&6 \end{array} \right] &
\left[ \begin{array}{rrrr}3&1&0&-1 \\ 1&3&1&0 \\ 0&1&6&2 \\ -1&0&2&6
\end{array}\right] & \mbox{---} \\ 
23 & 2 & \left[\begin{array}{cc}
4&1 \\ 1&6 \end{array}\right] & \left[\begin{array}{cc}
3 & 1 \\ 1 & 8 \end{array}\right] & \mbox{---} \\  \hline
\end{array}
$$
\end{table}
\begin{theorem}\label{ThCon1}
Consider the Mathieu group $M_{23}$ acting on the Leech lattice
$\La_{24}$, and let $g \in M_{23}$ have order $N > 1$.
(There is essentially only one class of elements of each order.)
Then the sublattice $\La_g$ of $\La_{24}$ fixed by $g$ is strongly $N$-modular.
If $N$ is in \eqn{Eq4} then $\La_g$ is extremal of dimension $D_N$.
\end{theorem}
\paragraph{Proof.}
A straightforward case-by-case verification.
(The $\La_g$ are also described in
\cite{HL90},
\cite{Koike85},
\cite{KT86}.)~~~$\bsq$
\paragraph{Remarks.}
(1)~We were led to this result by Quebbemann's observation in \cite{Q95} (following \cite{Con17}) that the function field of
$\Gamma_0 (p)^+$ for $p$ prime has genus 0 exactly when $p$ divides
the order of the Monster group.
Our investigations had suggested the group $\frac{1}{2} \Gamma_0 (4p)^+$ and the list of primes 2, 3, 5, 7, 11, 23.
It was natural to conjecture that these primes also arose from some finite simple group, the
obvious candidates being $M_{23}$, $M_{24}$ and $Co_2$.
The theta series of the sublattices $\La_g$ given by
Koike \cite{Koike85} then suggested the theorem.

(2)~It has been observed (\cite{HL89}, \cite{KT86}, \cite{KT87}) that the theta series of the fixed sublattices $\La_g$ for $g \in Co_0$ transform nicely
under Atkin-Lehner involutions.
For $g \in M_{23}$ this can be independently deduced from the modularity of $\La_g$,
using Corollary~\ref{Co3} (this does not seem to have been noticed before).
Indeed, it turns out that every relation between these theta series
under Atkin-Lehner involutions can be explained by an appropriate modularity.

(3)~There are two conjugacy classes of $M_{23}$ with orders not in \eqn{Eq4}, 
those of orders 4 and 8.
For order 4 the fixed sublattice is the 10-dimensional 4-modular
lattice called $Q_{10}$ in \cite{Me177}, \cite{HL90}.
For order 8 it is a 6-dimensional 8-modular even lattice with minimal
norm 4 and automorphism group of order 384 \cite{HL90}.

\begin{theorem}\label{ThCon2}
(a) Consider $M_{23}$ acting on the odd Leech lattice $O_{24}$,
and suppose $g \in M_{23}$ has odd order $N > 1$.
The fixed sublattice $\La_g$ is a strongly
$N$-modular $D_N$-dimensional lattice $O^{(N)}$ of minimal norm 3.
(b)~Consider $M_{23}$ acting on $O_{23} \oplus \ZZ$.
Again supposing that $N$ is odd, the fixed sublattice has the form
$\CN \oplus D$, where $D$ has dimension $0$ if $N=15$ or $23$,
and otherwise is an $N$-modular lattice $S^{(N)}$ of dimension $D_N - \dim \, \CN$ and minimal norm $3$.
\end{theorem}
\paragraph{Proof.}
Again a case-by-case verification.~~~$\bsq$

Since there is no exceptional case in Theorem~\ref{th2} when $N$
is even, the shorter lattices $S^{(N)}$ do not exist.
There are however odd lattices $O^{(N)}$ for $N=2$, 6 and 14,
although the construction of Theorem~\ref{ThCon2}
does not work.
The most interesting of these cases is $N=2$, for which
$S^{(2)}$ can be constructed as follows.

Let $L$ denote the 16-dimensional 2-modular lattice $BW_{16}$, with
minimal norm 4, and take $v \in L$ with $v \cdot v = 6$, $w \in L^\ast$ with $w \cdot w =3$.
Then
$O^{(2)} = \langle L', w \rangle$ where
$L' = \{ u \in L : u \cdot v \in 2 \ZZ \}$.

In fact all the $O^{(N)}$ and $S^{(N)}$ in Table~\ref{TaE}
can be found by a similar neighboring process, starting
from the even lattice $E^{(N)}$.
In each case there are four equivalence classes of $E^{(N)} / 2E^{(N)}$ under
the action of $\mathop{\rm Aut} (E^{(N)})$, with minimal norms 0, 4, 6, 8.
Relative to a vector of norm 4, the even neighbor is $E^{(N)}$
again, and the odd neighbor is $\CN \oplus S^{(N)}$.
Relative to a vector of norm 8,
the even neighbor is an analog of the Niemeier lattice of type
$A_1^{24}$, while the odd neighbor is $O^{(N)}$.

All the lattices in Table~\ref{TaE} are unique, although we only discuss
the uniqueness of $O^{(N)}$ and $S^{(N)}$ here.
It can be shown that if $\La$ is any $N$-modular lattice of norm 3 in the same
dimension as $S^{(N)}$, then the even neighbors of $\CN \oplus S^{(N)}$ must be extremal;
this implies the uniqueness of $S^{(N)}$.
A similar argument (based on the fact that the analogue of $A_1^{24}$ has the minimal nonzero number of roots) shows the uniqueness of $O^{(N)}$ for $N$ odd.
For $N$ even, one can show (see Theorem~\ref{ThEN}) that the even neighbor
of such a lattice must be extremal, and again the uniqueness of $O^{(N)}$ follows.

Finally, we comment on some of the other entries in Table~\ref{TaE}.
The 5-modular lattices $Q_8 (1)$ and $Q_6 (4)^{+2}$ are
connected with the ring of icosian integers --- see \cite{LDL2} (and \cite{PP8}).
$O^{(5)}$, $O^{(6)}$ and $O^{(7)}$ may be new:
they have minimal norm 3,
automorphism groups of orders 384, 96 and 48,
respectively, and 16, 16 and 8 minimal vectors
(see \cite{SN}).
The remaining entries are self-explanatory.

\section{Modular lattices and Atkin-Lehner involutions}
\hsp
A lattice $\La$ is {\em rational} (resp. {\em integral}) if $u \cdot v \in \QQ$
(resp. $\ZZ$) for all $u,v \in \La$.
Let $\Pi$ be a (possibly infinite) set of rational primes.
The {\em $\Pi$-dual} $\La^{\ast_\Pi}$ of $\La$ consists of
the vectors $v \in \La \otimes \QQ$ such that
$v \cdot \La \subseteq \ZZ_p$ for $p \in \Pi$ and $v \cdot \LA \subseteq \ZZ_p$ for $p \not\in \Pi$.

In particular, with $\Om$ the set of all rational primes,
$$
\La^{\ast_\emptyset} = \La , \quad
\La^{\ast_\Omega} = \LA , \quad
(\LAP )^{\ast_\Pi} = \La ~,
$$
and, more generally,
$$( \La^{\apo} )^{\apt} = \La^{\ast_{(\Pi_1\Delta \Pi_2)}}$$
where $\Delta$ denotes a symmetric difference.
(We will also need the notation $\overline{\Pi} = \Omega \setminus \Pi$,
and when there is no possibility of confusion we abbreviate
$\Pi = \{p\}$ to $p$.)
Furthermore,
$$\LAP \otimes \ZZ_p = \left\{
\begin{array}{ll}
\La^\ast \otimes \ZZ_p , & p \in \Pi ~, \\
\La \otimes \ZZ_p , & p \not\in \Pi ~.
\end{array}
\right.
$$
We also define
\begin{eqnarray*}
{\rm det}_\Pi (\La ) & = &
( {\rm det}~ \La / {\rm det}~ \LAP )^{1/2} \\
& = &
[\LAP : \La \cap \LAP ] / [ \La : \La \cap \LAP ] ~,
\end{eqnarray*}
which is equal to the $\Pi$-part of ${\rm det}~ \La$.

Suppose now that $\La$ is integral.
The {\em level} of $\La$ is the smallest number $l'$ such that $\sqrt{l'} \La^\ast$ is integral.
If $\La$ is even, the {\em even-level} of $\La$ is the smallest
number $l$ such that $\sqrt{l} \La^\ast$ is even.
The $\Pi$-levels $l'_\Pi$ and $l_\Pi$ are defined
analogously, replacing $\LA$ by $\LAP$.

Quebbemann \cite{Q97} associates certain Gauss sums with $\La$.
We do the same, but in a slightly more explicit fashion.
Let
$$\gamma_2 (\La ) = \xi^{{\rm oddity} ( \La )} , \quad
\gamma_p ( \La ) = \xi^{-{p{\rm \hbox{-}excess}} ( \La )} ~,
$$
for an odd prime $p$, where $\xi = e^{\pi i/4}$ and the oddity and $p$-excess
are as in Chap. 15 of \cite{SPLAG}, and define
\beql{Eq13a}
\gamma_\Pi ( \La ) = \prod_{p \in \Pi} \gamma_p ( \La ) ~.
\eeq
In particular, the product (or oddity) formula
\cite[Chap. 15, Eq. (30)]{SPLAG} becomes
$$\gamma_\Om ( \La ) = \xi^{\dim \La} ~.$$
The following lemma shows that $\gamma_\Pi ( \La )$ agrees with Quebbemann's 
Gauss sum.
\begin{lemma}\label{Lem1}
For an even lattice $\La$,
\beql{Eq13b}
\gamma_\Pi ( \La ) = ({\rm det}_\Pi~
\La)^{-1/2} \sum_{v \in \LAP/ \La} e^{\pi i v \cdot v} ~.
\eeq
\end{lemma}
\paragraph{Proof.}
From \cite[Chap. 5]{Sch85}, the right-hand side of \eqn{Eq13b}
is multiplicative under direct sums of lattices and disjoint
unions of prime sets,
and is invariant under rational equivalence of lattices.
It suffices therefore to consider only the cases where $\Pi$
is a singleton and $\La = \sqrt{a} \ZZ$, where $a$ ranges over
$\ZZ_p^\ast / (\ZZ_p^\ast )^2$.
This is a straightforward problem involving one-dimensional Gauss sums.~~~$\bsq$

It is classical (cf. \cite{Mi}) that if $\La$ is a lattice of
even-level $N$, then its theta series $\Theta_\La$ is a modular form
for $\Gamma_0 (N)$ with respect to an appropriate character.  Kitaoka
\cite{Kit80} describes how a somewhat larger subset of $SL_2 ( \ZZ )$
acts on $\Theta_\La$, up to an unspecified constant.  Quebbemann
\cite{Q97} has determined this constant, but only for one
representative from each coset of $\Gamma_0 (N)$.  We shall make use
of the following more explicit result.  Here $\Pi (m)$ denotes the set
of primes dividing $m$, and $\legendre{m}{n}$ denotes the
Kronecker-Jacobi symbol \cite[p. 28]{Cohen}.

\begin{theorem}\label{ThAL}
Let $\Lambda$ be an even lattice of even-level $N$, and let
$S=\left(\begin{array}{cc}a & b \\ c & d \end{array}\right)$ be any element of $SL_2({\Bbb Z})$ such that $cd$ is a
multiple of $N$.  Then
\beql{Eq13c}
\Theta_\Lambda((az+b)/(cz+d)) =
({\rm det}_{\Pi(d)}~\Lambda)^{-1/2}
\chi_{c,d}(\Lambda)
\left(\sqrt{cz+d}\right)^{\dim\Lambda}
\Theta_{\Lambda^{*_{\Pi(d)}}}(z),
\eeq
where in both cases the square root is that with positive real part, and
$\chi_{c,d}(\Lambda)$ is equal to
$$
\gamma_{\Pi(d)}(\Lambda)^{-1}
\legendre{c}{{\rm det}_{\Pi(d)}~\Lambda}
\legendre{d}{{\rm det}_{\Pi(c)}~\Lambda}
$$
multiplied either by
\beql{Eq13d}
\legendre{d}{|c|}^{\dim\Lambda}
\legendre{\legendre{-1}{c}}{{\rm det}~\Lambda} \xi^{-(c-1)\dim\Lambda}
\eeq
if $c$ is odd, or by
\beql{Eq13e}
\legendre{c}{d}^{\dim\Lambda}
\legendre{\legendre{-1}{d}}{{\rm det}~\Lambda} \xi^{(d-1)\dim\Lambda}.
\eeq
if $c$ is even.
\end{theorem}

For the proof, we need a
lemma describing how Gauss sums behave as a lattice is
rescaled.

\begin{lemma}\label{LemAL}
Let $\Lambda$ be a rational lattice, and let $\Pi$ be any set of primes.
Let $t$ be any positive integer, with $\Pi$-part $t_1$ and
$\overline{\Pi}$-part $t_2$.  If $2\not\in \Pi$, then
\beql{Eq13f}
\gamma_\Pi(\sqrt{t} \Lambda)/\gamma_\Pi(\Lambda)
=
\legendre{t_1}{{\rm det}_{\overline{\Pi}}~\Lambda}
\legendre{t_2}{{\rm det}_\Pi~\Lambda}
\legendre{t_2}{t_1}^{\dim\Lambda}
\legendre{\legendre{-1}{t_1}}{{\rm det}~\Lambda} \xi^{-(t_1-1)\dim\Lambda},
\eeq
and if $2\in \Pi$, then
\beql{Eq13g}
\gamma_\Pi(\sqrt{t} \Lambda)/\gamma_\Pi(\Lambda)
=
\legendre{t_1}{{\rm det}_{\overline{\Pi}}~\Lambda}
\legendre{t_2}{{\rm det}_\Pi~\Lambda}
\legendre{t_1}{t_2}^{\dim\Lambda}
\legendre{\legendre{-1}{t_2}}{{\rm det}~\Lambda} \xi^{(t_2-1)\dim\Lambda}.
\eeq
\end{lemma}

\paragraph{Proof.}
It follows from the definition of the $p$-excess that if $p$
is any odd prime and $t$ is relatively prime to $p$ then
$$
\gamma_p(\sqrt{t} \Lambda)/\gamma_p(\Lambda)
=
\legendre{t}{{\rm det}_p~\Lambda}.
$$
Furthermore, if $p$ does not divide ${\rm det}~\Lambda$, then
\begin{eqnarray*}
\gamma_p(\sqrt{p} \Lambda)/\gamma_p(\Lambda) &=&
\xi^{-(p-1)\dim\Lambda}
\legendre{{\rm det}~\Lambda}{p}\\
&=&
\xi^{-(p-1)\dim\Lambda} \legendre{\legendre{-1}{p}}{{\rm det}~\Lambda}
\legendre{p}{{\rm det}~\Lambda},
\end{eqnarray*}
by reciprocity.

For $p=2$ and $t$ odd, the result clearly depends only on the congruence
class of $t$ mod 8.  Consequently, we may assume that $t$ is a prime not
dividing ${\rm det}~\Lambda$.  Then
\begin{eqnarray*}
\gamma_2(\sqrt{t} \Lambda)/\gamma_2(\Lambda)
&=&
(\gamma_t(\sqrt{t} \Lambda)/\gamma_t(\Lambda))^{-1}
(\gamma_{\overline{\{2,t\}}}(\sqrt{t} \Lambda)/
 \gamma_{\overline{\{2,t\}}}(\Lambda))^{-1}\\
&=&
\xi^{(t-1)\dim\Lambda} \legendre{\legendre{-1}{t}}{{\rm det}~\Lambda}
\legendre{t}{{\rm det}~\Lambda}\legendre{t}{{\rm det}_{\overline{\{2,t\}}}~\Lambda}\\
&=&
\xi^{(t-1)\dim\Lambda} \legendre{\legendre{-1}{t}}{{\rm det}~\Lambda}
\legendre{t}{{\rm det}_2~\Lambda}.
\end{eqnarray*}

We can now write, for $2\not\in\Pi$:
$$
\gamma_\Pi(\sqrt{t} \Lambda)/\gamma_\Pi(\Lambda)
=
(\gamma_\Pi(\sqrt{t_1t_2} \Lambda)/\gamma_\Pi(\sqrt{t_1}\Lambda))
(\gamma_\Pi(\sqrt{t_1} \Lambda)/\gamma_\Pi(\Lambda)).
$$
The first ratio is
$$
\legendre{t_2}{{\rm det}_\Pi (\sqrt{t_1} \Lambda)}
=
\legendre{t_2}{t_1}^{\dim\Lambda} \legendre{t_2}{{\rm det}_\Pi~\Lambda}
$$
while the second is
\begin{eqnarray*}
\gamma_\Pi(\sqrt{t_1} \Lambda)/\gamma_\Pi(\Lambda)
&=&
(\gamma_{\overline{\Pi}}(\sqrt{t_1} \Lambda)/\gamma_{\overline{\Pi}}(\Lambda))^{-1}\\
&=&
(\gamma_2(\sqrt{t_1} \Lambda)/\gamma_2(\Lambda))^{-1}
(\gamma_{\overline{\Pi\cup\{2\}}}(\sqrt{t_1} \Lambda)
/\gamma_{\overline{\Pi\cup\{2\}}}(\Lambda))^{-1}\\
&=&
\xi^{-(t_1-1)\dim\Lambda} \legendre{\legendre{-1}{t_1}}{{\rm det}~\Lambda}
\legendre{t_1}{{\rm det}_2~\Lambda}
\legendre{t_1}{{\rm det}_{\overline{\Pi\cup\{2\}}}~\Lambda}\\
&=&
\xi^{-(t_1-1)\dim\Lambda} \legendre{\legendre{-1}{t_1}}{{\rm det}~\Lambda}
\legendre{t_1}{{\rm det}_{\overline{\Pi}}~\Lambda}.
\end{eqnarray*}
This establishes \eqn{Eq13f}.  \eqn{Eq13g} then follows from the oddity
formula.~~~$\bsq$

\paragraph{Proof of Theorem~\ref{ThAL}.}
We first suppose $c>0$.  Quebbemann
\cite{Q97} shows that when $a=1$ and $c|N$,
$$
\Theta_\Lambda((z+b)/(cz+d)) =
({\rm det}_{\Pi(d)}~\Lambda)^{-1/2}
\xi^{-\dim\Lambda} \gamma_{\Pi(c)}(\sqrt{c}\Lambda)
\left(\sqrt{cz+d}\right)^{\dim\Lambda}
\Theta_{\Lambda^{*_{\Pi(d)}}}(z).
$$
(To be precise, \cite{Q97} has $\gamma_{\Pi(c)}(\sqrt{c}\Lambda^{*_{\Pi(c)}})$,
but since $\Lambda$ and $\Lambda^{*_{\Pi(c)}}$ are rationally equivalent,
this is the same as $\gamma_{\Pi(c)}(\sqrt{c}\Lambda)$.)
The argument in \cite{Q97} never uses the fact that $c$ divides $N$, and can be
easily modified to show that for arbitrary $a>0$,
$$
\Theta_\Lambda((az+b)/(cz+d)) =
({\rm det}_{\Pi(d)}~\Lambda)^{-1/2}
\xi^{-\dim\Lambda} \gamma_{\Pi(c)}(\sqrt{ac}\Lambda)
\left(\sqrt{cz+d}\right)^{\dim\Lambda}
\Theta_{\Lambda^{*_{\Pi(d)}}}(z).
$$
If $c$ is odd, the lemma implies
\begin{eqnarray*}
\xi^{-\dim\Lambda} \gamma_{\Pi(c)}(\sqrt{ac} \Lambda)
&=&
\left(\xi^{-\dim\Lambda} \gamma_{\Pi(c)}(\Lambda)\right)
\legendre{c}{{\rm det}_{\overline{\Pi(c)}}~\Lambda}
\legendre{a}{{\rm det}_{\Pi(c)}~\Lambda} \cdot\\
&&~~~~~~~~~~~\cdot~
\legendre{a}{c}^{\dim\Lambda}
\legendre{\legendre{-1}{c}}{{\rm det}~\Lambda} \xi^{-(c-1)\dim\Lambda} \cdot\\
&=& \gamma_{\Pi(d)}(\Lambda)^{-1}
\legendre{c}{{\rm det}_{\Pi(d)}~\Lambda}
\legendre{d}{{\rm det}_{\Pi(c)}~\Lambda} \cdot\\
&&~~~~~~~~~~\cdot~
\legendre{d}{c}^{\dim\Lambda}
\legendre{\legendre{-1}{c}}{{\rm det}~\Lambda} \xi^{-(c-1)\dim\Lambda},
\end{eqnarray*}
where the second step follows from the oddity formula and the fact that
$ad\bmod c=1$.  For $a \leq 0$, we use the fact that
$\Lambda$ is even, so the result can depend only on the value of $a$ mod
$c$.

For $c$ even, we do not, in general have $\legendre{a}{2}=\legendre{d}{2}$,
so the above argument fails.  However, again using the fact that the result
only depends on the value of $a$ mod $c$, we can arrange that $a\equiv
d\pmod{8}$, and then an analogous argument can be used.

For $c$ negative ($c=0$ is trivial), we apply the result to $-S$,
and use the fact that $\sqrt{-cz-d}=i\sqrt{cz+d}$.
For $c$ odd,
$$
i^{\dim\Lambda} \chi_{-c,-d}(\Lambda)
=
\chi_{c,d}(\Lambda)
\left(
-i \legendre{-1}{c} \xi^{2c}\right)^{\dim\Lambda},
$$
while for $c$ even,
$$
i^{\dim\Lambda} \chi_{-c,-d}(\Lambda)
=
\chi_{c,d}(\Lambda)
\left(
i \legendre{-1}{d} \xi^{-2d}\right)^{\dim\Lambda}.
$$
But for odd integers $n$, $\legendre{-1}{n} \xi^{2n}=i$, so in either
case,
$$
\chi_{c,d}(\Lambda)\left(\sqrt{cz+d}\right)^{\dim\Lambda}=
\chi_{-c,-d}(\Lambda)\left(\sqrt{-cz-d}\right)^{\dim\Lambda},
$$
and so the above formulae also hold if $c$ is negative.~~~$\bsq$

\paragraph{Remarks.}
(1) There is an apparent inconsistency in \eqn{Eq13c}.  Since
$$
\Theta_{\Lambda^{*_{\Pi(d)}}}\left(z+{N\over \gcd(c,N)}\right)
=
\Theta_{\Lambda^{*_{\Pi(d)}}}(z),
$$
$\chi_{c,d}(\Lambda)$ must be periodic in $d$ of period $cN/\gcd(c,N)$.
For $c$ odd or $c\equiv 0\pmod{8}$ this is manifestly true, but otherwise
\eqn{Eq13e} appears to have the wrong period.  For instance, for $c\equiv
4 ~ (\bmod ~ {8})$,
$$
\chi_{c,d+(cN/\gcd(c,N))}(\Lambda) =
(-1)^{\lambda+\dim\Lambda}~\chi_{c,d}(\Lambda),
$$
where $\lambda=\log_2({\rm det}_2~\Lambda)$.  However, since $N|cd$, it follows
that in the $2$-adic Jordan decomposition of $\Lambda$ the forms of levels
$1$ and $4$ are both Type II and so have even dimension.  This implies that
$\lambda\equiv\dim\Lambda~ (\bmod ~ {2})$.

Similarly, for $c\equiv{\pm 2} ~ (\bmod ~{8})$, the correct period is restored by
the identities
$$
\displaylines{
\lambda\equiv\dim\Lambda\equiv 0\pmod{2},\cr
\legendre{-1}{{\rm det}~\Lambda}=(-1)^{(\dim\Lambda)/2}.\cr}
$$
If both $c$ and $d$ are odd (so
$\Lambda$ has odd even-level), then similar reasoning allows us to simplify
$\chi_{c,d}(\Lambda)$ to
$$
\chi_{c,d}(\Lambda)
=
\gamma_{\Pi(d)}(\Lambda)^{-1}
\legendre{c}{{\rm det}_{\Pi(d)}~\Lambda}
\legendre{d}{{\rm det}_{\Pi(c)}~\Lambda}.
$$

(2) When $N$ divides $c$, the usual formula \cite[Theorem 4.9.3]{Mi}
for the action of $\Gamma_0(N)$ on the theta series of lattices of
even-level $N$ can be recovered with the help of the identity
$$
\xi^{(t-1)}=\epsilon_t \legendre{-2}{t}=\epsilon_t^{-1} \legendre{2}{t},
$$
for odd $t$, where $\epsilon_t=1$ if $t\equiv 1\pmod{4}$ and $\epsilon_t=i$
if $t\equiv 3\pmod{4}$.

If $\Lambda$ is any integral lattice of level $N$, $\sqrt{2}\Lambda$
is an even lattice of even-level dividing $4N$.  We can apply Theorem
\ref{ThAL} to obtain:

\begin{coro}\label{CoAL}
Let $\Lambda$ be an integral lattice of level $N$, and let
$S=\left( \begin{array}{cc} a&b\\c&d\end{array}\right)$ be any element of $SL_2({\Bbb Z})$ such
that $cd$ is a multiple of $2N$.  Then $($\ref{Eq13c}$)$ holds if either $d$ is odd
and $b$ is even, or $c$ is odd and $a$ is even.
\end{coro}

A {\em modularity} $\sigma$ of an integral lattice $\Lambda$ is a similarity
mapping $\Lambda^{*_\Pi}$ to $\Lambda$ for some
set of primes $\Pi$.  We say that
$\sigma$
has {\it level} $N$ (or is an {\it $N$-modularity}) if $\sigma$ multiplies
norms by $N$; $\Pi$ is then the set of primes dividing $N$.  A
$1$-modularity is just an automorphism of $\Lambda$.

\begin{coro}\label{Co3}
Suppose $\Lambda$ has even-level $N$ and admits an
$m$-modularity.  Then for any matrix
\beql{Eq13h}
W_m=m^{-1/2} \left(\begin{array}{cc} ma&b\\mc&d\end{array}\right)
\eeq
of determinant $1$, with $d$ a multiple of $m$ and $mc$ a multiple of $N$,
we have
$$
\Theta_\Lambda|_{W_m} = \chi_{c,d}(\Lambda) \Theta_\Lambda.
$$
\end{coro}

\paragraph{Proof.}
Note that
$$
{\rm det}_{\Pi(m)}~\Lambda
= \left({{\rm det}~\Lambda\over {\rm det}~\Lambda^{*_{\Pi(m)}}}\right)^{1/2}
= {\rm det}~\sigma = m^{(\dim\Lambda)/2},
$$
and $\sqrt{m}\Lambda^{*_{\Pi(m)}}$ is isometric to $\Lambda$, where the
isometry is $\sigma/\sqrt{m}$.  Applying Theorem \ref{ThAL}, we find
\begin{eqnarray*}
\Theta_\Lambda((amz+b)/(cmz+d)) &=&
({\rm det}_{\Pi(d)} \Lambda)^{-1/2}
\chi_{c,d}(\Lambda)
\left(\sqrt{cmz+d}\right)^{\dim\Lambda}
\Theta_{\Lambda^{*_{\Pi(d)}}}(mz)\\
&=&
m^{-(\dim\Lambda)/4}
\chi_{c,d}(\Lambda)
\left(\sqrt{cmz+d}\right)^{\dim\Lambda}
\Theta_{\sqrt{m} \Lambda^{*_{\Pi(m)}}}(z)\\
&=&
\chi_{c,d}(\Lambda)
\left(\sqrt{m^{1/2}cz+m^{-1/2}d}\right)^{\dim\Lambda}
\Theta_{\Lambda}(z). ~~~~~\bsq
\end{eqnarray*}

The matrix $W_m$ in \eqn{Eq13h} is called an
{\em Atkin-Lehner involution} \cite{AL} of level $m$.
The next result combines known properties of these involutions
with a slight generalization of a result of Nebe \cite{NebeSM} on modularities.
We omit the proof.
\begin{theorem}\label{ThPm}
If $W_{m_1}$ and $W_{m_2}$ are Atkin-Lehner involutions then $W_{m_1} W_{m_2}$ is an Atkin-Lehner involution of level $m_1 m_2 / gcd (m_1, m_2)^2$.
Moreover, $W_m^{-1}$ is an Atkin-Lehner involution of level $m$.
If $\sg_1$ is an $m_1$-modularity and $\sg_2$ is an $m_2$-modularity
then $\sg_1 \sg_2 / gcd (m_1, m_2)$ is a modularity of level
$m_1 m_2 / gcd (m_1, m_2)^2$.
Moreover, if $\sg$ is an $m$-modularity then so is $m\sg^{-1}$.
\end{theorem}

It follows from Theorem~\ref{ThPm} that the number of distinct levels of modularities of a lattice is a power of 2, and indeed the levels have a natural elementary abelian
2-group structure.
Moreover, the total number of modularities is equal to the number of levels
of modularity times $|\mathop{\rm Aut}\La |$.

We will say that an integral lattice $\La$ is
{\em $\{l_1, l_2, \ldots \}$-modular} if it has modularities of levels
$l_1 , l_2, \ldots$.
Two special cases warrant a shorthand notation.
(i)~$\La$ is {\em $N$-modular} if its level divides $N$ and $\La$ is
$\{1,N\}$-modular.
(ii)~$\La$ is {\em strongly $N$-modular} if its level divides $N$ and $\La$
is $\{m: m \| N \}$-modular, where $a \| b$ means $a|b$ and $gcd (a, b/a) =1$.

Corollary~\ref{Co3} states that if $\La$ is an even $\{l_1, l_2, \ldots \}$-modular lattice of even-level $N$,
then its theta series is an automorphic form for the
group $\Gamma_0 (N)^{+ \{ l_1, l_2 , \ldots \}}$, i.e. the group
generated by $\Ga_0 (N)$ together with all its Atkin-Lehner involutions
of levels $l_1 , l_2 , \ldots$.
For ease in discussing strongly modular lattices we abbreviate
$\Ga_0 (N)^{+ \{m: m \| N \}}$ to $\Ga_0 (N)^+$.

If $\Ga$ is any modular group, $\frac{1}{2} \Ga$ will denote the
group $\left \{ \left( \begin{array}{cc}
a & 2b \\ c/2 & d
\end{array} \right) : \left( \begin{array}{cc}
a & b \\ c & d \end{array} \right) \in \Ga \right\}$.

Corollary \ref{CoAL} implies that if $\La$ is an
$\{l_1, l_2, \ldots \}$-modular lattice of level $N$, then its theta series
is an automorphic form for the group
\beql{Eq23a}
\frac{1}{2} \Ga_0 (4N)^{+ \{ 4, l_1, l_2, \ldots \}} , \quad
\mbox{if $N$ is odd} ~.
\eeq
The initial 4 arises because $\sqrt{2} \La$ has an obvious 4-modularity.
As a special case, the theta series of any lattice of odd level is an
automorphic form for $\frac{1}{2} \Ga_0 (4N)^{+ \{4\}}$ (a subgroup of $\Ga_0 (N)$).
If $N$ is even, \eqn{Eq23a} must be replaced by
$$\frac{1}{2} \Gamma_0 (4N)^{+ \{4e_1, 4e_2 , \ldots, d_1, d_2, \ldots \}} ~,$$
where $e_1, e_2 , \ldots$ are the even $l_i$'s and $d_1$, $d_2, \ldots$ are the odd $l_i$'s.

\section{Shadows}
\hsp
Let $\La$ be an integral lattice, or more generally a 2-integral lattice (i.e. $u \cdot v \in \ZZ_2$ for all $u,v \in \La$), and set $\La_0 = \{u \in \La : u \cdot u \in 2 \ZZ_2 \}$.
If $\La$ is even, $\La = \La_0$;
otherwise $\La_0$ is a sublattice of index 2.
$\La_0$ is called the even sublattice of $\La$.

Following \cite{Me157}, \cite{Me158}, we define the {\em shadow}
$S( \La )$ to be $\La^\ast$ if $\La$ is even,
$(\La_0)^\ast \setminus \LA$ if $\La$ is odd.
Equivalently,
$$S( \La ) = \{ v \in \La \otimes \QQ ~:~
2 u \cdot v \equiv u \cdot u~ ( \bmod~2 \ZZ_2 ) \quad\mbox{for all}\quad
u \in \La \} ~.
$$
Also
\beql{Eq14a}
\Theta_{S( \La )} (z) =
({\rm det} ~ \La)^{1/2} \left(
\frac{\xi}{\sqrt{z}} \right)^{\dim \La} \Theta_{\La} \left( 1- \frac{1}{z} \right) ~.
\eeq 

We also define the {\em $\Pi$-shadow} $S_\Pi (\La)$: if $2 \in \Pi$,
$$S_{\Pi} ( \La ) = S(\La^{\ast_{\overline{\Pi}}} ) ~,$$
and if $2 \not\in \Pi$,
$$S_\Pi ( \La ) = \sqrt{l'_2} S( \sqrt{l'_2} \La^{\ast_{\overline{\Pi}}} )$$
where $l'_2$ is the 2-level of $\La$.
The $\Pi$-shadow is a coset of the $\Pi$-dual $\LAP$, and in fact
$v \pm w \in \LAP$ for $v,w \in S_\Pi ( \La )$.
In particular,
$S_\Omega (\La) = S( \La )$ is a coset of $\LA$, and $S_\emptyset ( \La )$ is a coset of $\La$.
The theta-series of $S_\Pi ( \La )$ may be computed from Corollary~\ref{CoAL} and \eqn{Eq14a}.

It is clear from the definition of $S(\La)$ that any two vectors
in the same coset of $\La$ in $S( \La)$ have the same norm modulo $2 \ZZ_2$.
If $\La$ has odd determinant, we can say more.
\begin{theorem}\label{ThOdd}
Let $\La$ be a $2$-integral lattice of odd determinant and let $\Pi$ be a set of rational primes.
Then every vector in $S_\Pi (\La )$ has $\mbox{norm $\equiv$ $($oddity $\La)/4$}$
$(\bmod~2 \ZZ_2 )$.
\end{theorem}
\paragraph{Proof.}
We give three proofs.
It suffices to consider $\Pi = \Omega$, since
$\La^{\ast_{\overline{\Pi}}}$ satisfies the hypotheses and has the same oddity.

First proof.
By scaling $\La$ we may assume $\La$ is integral.
Since $\La$ has odd determinant, $\La^{\ast_2} = \La$.
Applying Corollary~\ref{CoAL}, we have
$$\Theta_\La \Bigl|_{\left[ \begin{array}{cc}
1+N^2 & -N^2 \\ N^2 & 1-N^2
\end{array} \right]} =
\chi_{N^2,1-N^2} ( \La ) \Theta_{\La^{\ast_2}} ~.
$$
Now $\chi_{N^2, 1-N^2} (\La ) = \gamma_2 ( \La)^{-1}$.
Since
$$\left[ \begin{array}{cc}
1+N^2 & -N^2 \\ N^2 & 1-N^2
\end{array}
\right] ~~
\left[ \begin{array}{cr}
1 & -1 \\ 1 & 0
\end{array}
\right] = \left[ \begin{array}{cr}
1 & -1 \\ 1 & 0
\end{array}
\right] ~~
\left[ \begin{array}{cc}
1 & -N^2 \\ 0 & 1
\end{array}\right] ~,
$$
we have
$$\Theta_{S(\La)} (z-N^2) = \gamma_2 (\La)^{-1} \Theta_{S( \La )} (z) ~.
$$
In other words,
$$e^{- \pi i N^2 v \cdot v} = e^{- \frac{2 \pi i}{8} ~{\rm oddity \La}}$$
for all $v \in S( \La )$.

Second proof.
Since the desired result is purely 2-adic,
we may localize at the prime 2.
Because $S( \La_1 \oplus \La_2 ) = S(\La_1 ) \oplus S(\La_2)$, the
result is preserved under direct summation, so it suffices
to consider indecomposable 2-adic quadratic forms.
It is straightforward to verify that the theorem holds for each
of the six classes of 1- or 2-dimensional forms of unit determinant.

Third proof.
Assume $\La$ is integral and odd ($\La$ even is trivial).
Since $\La_0 \subseteq \La$ and oddity is a rational invariant,
\begin{eqnarray*}
\gamma_2 ( \La ) = \gamma_2 ( \La_0) & = &
( {\rm det}_2 ~ \La_0 )^{-1/2} \sum_{v \in \La_0^{\ast_2} / \La_0}
e^{\pi i v \cdot v} ~ \\
& = & \frac{1}{2}
\sum_{v \in \La / \La_0} e^{\pi i v \cdot v} + 
\sum_{v \in S_2( \La ) / \La} e^{\pi i v \cdot v} ~.
\end{eqnarray*}
The first sum is $1-1 =0$, so
$e^{\pi i v \cdot v} = e^{\pi i~{\rm oddity}/4}$.~~~$\bsq$

\paragraph{Remarks.}
(1) For unimodular lattices, Theorem~\ref{ThOdd} together with the product
formula implies that for $v\in S(\Lambda)$, $v\cdot v\equiv
(\dim\Lambda)/4\pmod{2}$,
a result that has been rediscovered several times
(see \cite{Blij59}, \cite{Brau40}, \cite{Me157}, \cite{MH73}, \cite{Serre}).
(2) The third proof can be used to extend Lemma~\ref{Lem1} to integral
lattices, since it proves that
$$\gamma_2 (\La) = ({\rm det}_2 ~ \La)^{-1/2} \sum_{v \in S_2( \La )/\La}
e^{\pi i v \cdot v} ~.
$$

\paragraph{Genus of $\La_0$.}
Assume $\La$ is odd.
Since the even sublattice $\La_0$ is defined 2-adically, its genus can be computed from that of $\La$.
(There is no change in the $p$-adic genus for $p \neq 2$.)
Indeed, the change in the genus depends only on the unit form in the 2-adic
Jordan decomposition of $\La$.

When the oddity is not zero, the existence conditions for 1- and
2-dimensional forms \cite[Theorem 11 of Chap. 15]{SPLAG} and the fact
that oddity is a rational invariant leave just one possibility.  When
the oddity is zero, $\La_0$ has a form at level 2, which from the
existence conditions could be either Type I or II.  But by
Theorem~\ref{ThOdd}, every vector in $\La_0^\ast$ has integral norm.
It follows that the form at level 2 must be Type II.  We thus obtain
the list of transforms shown in Table~\ref{Tat} (using the notation of
\cite[Chap. 15]{SPLAG}).

To avoid undue proliferation of parentheses we adopt the conventions
that the operation $\Lambda\to\Lambda_0$ takes precedence over
$\Lambda\to\Lambda^{*_2}$, and both take precedence over
$\Lambda\to\sqrt{2}\Lambda$.  Thus $\sqrt{2}\Lambda_0^{*_2}$ means
$\sqrt{2}((\Lambda_0)^{*_2})$.

\begin{table}[htb]
\caption{Genus of $\La_0$ in terms of genus of $\La$.}
\label{Tat}
$$
\begin{array}{ll}
\mbox{genus ($\La$)} & \mbox{genus ($\La_0$)} \\ [+.1in]
\mbox{[$1^{\pm n}$]$_0$} & 1^{\pm (n-2)}: 2^2 \\ 
\mbox{[$1^{\pm n}$]$_1$} & 1^{\pm (n-1)}: [4^1]_1 \\
\mbox{[$1^{\pm n}$]$_2$} & 1^{\pm (n-2)} [2^2]_2 \\
\mbox{[$1^{\pm n}$]$_3$} & 1^{\mp (n-1)}: [4^{-1} ]_3 \\
\mbox{[$1^{\pm n}$]$_4$} & 1^{\mp (n-2)} : 2^{-2} \\
\mbox{[$1^{\pm n}$]$_5$} & 1^{\mp (n-1)} : [4^{-1} ]_5 \\
\mbox{[$1^{\pm n}$]$_6$} & 1^{\pm (n-2)} [2^2]_6 \\
\mbox{[$1^{\pm n}$]$_7$} & 1^{\pm (n-1)}: [4^1]_7
\end{array}
$$
\end{table}

\begin{theorem}\label{ThEN}
Let $\La$ be an odd $\{2\}$-modular lattice with dimension and oddity $o$ both divisible by $4$.
Then $\La' = \sqrt{2} \La_0^{\ast_2}$ is an integral lattice, and
$\La ''=(\La')'$ is an even $\{2\}$-modular lattice, rationally equivalent to $\La$.
In fact, every modularity of $\La$ is a modularity of $\La''$.
\end{theorem}

We call $\La ''$ the {\em even neighbor} of $\La$.
\paragraph{Proof.}
The 2-adic genus of $\La$ must be $[1^{n/2}~ 2^{n/2} ]_o$.
From Table~\ref{Tat},
the genera of $\La_0$ and $\La '$ are respectively
$1^{(n-4)/2}~ [2^{(n+4)/2} ]_o$ and
$[1^{(n+4)/2} ]_o~ 2^{(n-4)/2}$, and so $\La'$ is integral.
Then $(\La')_0$ has 2-adic genus $1^{n/2} : 2^{n/2}$ if $o=0$ and $1^{-n/2} : 2^{-n/2}$ if $o= 4$.
So $\La''$ is even.

If $\sg$ is a modularity of $\La$ of odd level, then it is still a modularity at each step of the construction.
If $\sg$ is a 2-modularity, then
$\La'' = \sg ( \sg \La_0^{\ast_2} )_0^{\ast_2}$.
But then $\La_0 \subseteq \La'' \subseteq \sg \La_0^{\ast_2}$,
and $\La_0 \subseteq \sg (\La'')^{\ast_2} \subseteq \sg \La_0^{\ast_2}$.
Since there is only one even lattice between $\La_0$ and $\sg \La_0^{\ast_2}$, and both $\La''$ and $\sigma (\La'')^{\ast_2}$ are even (from the genus of $\La''$),
it follows that they are the same lattice, and thus $\sg$ is a 2-modularity of $\La''$. 

The remaining modularities carry over to $\La'$ by Theorem \ref{ThPm}.
Since $\La \cap \La'' = \La_0$, $\La$ and $\La''$ are 
clearly rationally equivalent.~~~$\bsq$

The theta series of $\La''$ is
$$\frac{1}{2} \{ \Theta_\La (z) + \Theta_\La (z+1) + \Theta_{S_\emptyset(\La)} (z) +
\Theta_{S_\emptyset (\La)} (z+1) \} ~.
$$
In particular, if $\La$ is a 16-dimensional 2-modular lattice of minimal norm 3,
$\La''$ has minimal norm 4 and so (by \cite{Q95}) must be the Barnes-Wall lattice.
This forces the construction for the odd Barnes-Wall lattice given in Section 2.

\section{Theta series of strongly modular lattices}
\hsp
Throughout this section we assume that $\La$ is a strongly $N$-modular lattice for $N$ in \eqn{Eq4}.
As remarked in Section 3, if $\La$ is even then $\Theta_\La (z)$ is invariant under
$\Gamma_0 (N)^+$ with respect to a certain character depending only on the
rational equivalence class of $\La$.
In all cases $\Theta_\La (z)$ is invariant under $\frac{1}{2} \Gamma_0 (4N)^+$, again with
respect to some character.
In order to prove Theorems \ref{th1} and \ref{th2}
it is necessary to study
the space of modular forms for $\frac{1}{2} \Gamma_0 (4N)^+$.

Let $\chi^{(N)}$ be the character of $\frac{1}{2} \Gamma_0 (4N)^+$ with respect to which $\Theta_{\CN} (z)$ is invariant, and let $w^{(N)}$ be the weight of
$\Theta_{\CN} (z)$.
Then a lattice satisfying the hypotheses of Theorem \ref{th2} has theta series in the space
$$\sM_{kw^{(N)}} \left( \frac{1}{2} \Gamma_0 (4N)^+ , (\chi^{(N)} )^k \right)~.
$$
Of course $\Theta_{\CN} (z)^k$ is in this space.

In the sequel, we define the divisor of a modular form $f(z)$ in
${\cal M}_k(\Gamma,\chi)$ to be
$$
{1\over n} {\rm div}(f(z)^n),
$$
where $\chi^n=1$, and the divisor of a form with trivial character is defined
as in \cite[p. 51]{Mi}.

\begin{lemma}\label{P1}
For any square-free $N$, the divisor of $\Theta_{\CN} (z)$
with respect to $\frac{1}{2} \Gamma_0 (4N)^+$
is
$$
\frac{1}{8} \sum_{d|N} d \cdot {\bf 1} , ~~~\mbox{if $N$ odd} ~,\quad
\frac{1}{6} \sum_{d|N} d \cdot {\bf 1} , ~~~\mbox{if $N$ even} ~.
$$
\end{lemma}
\paragraph{Proof.}
The modular form
$$
\eta (z) = q^{1/12} \prod_{m=1}^{\infty} (1-q^{2m}), ~~~ q = e^{\pi i z },
$$ is zero only at $\QQ ~ \cup ~ \{ \infty \}$.  It follows that any
product or quotient of functions $\eta (az+b)$ for rational $a$ and
$b$ has no zeros or poles outside $\QQ ~ \cup ~ \{ \infty \}$.  In
particular, since $\Theta_{\ZZ}(z)=\theta_3 (z) = \eta (z)^5 / ( \eta
(z/2) \eta (2z))^2$, the same is true for $\Theta_{\CN} (z)$
($\theta_2(z)$, $\theta_3(z)$, and $\theta_4(z)$ are the familiar
Jacobi theta series).  Since $\CN$ is a lattice, $\Theta_{\CN} (z)$
does not have a zero at $\infty$.  Consequently,
$${\rm div}(\Theta_{\CN} (z)) = {\rm deg} (\Theta_{\CN} (z)) \cdot {\bf 1} ~.$$
We may compute the right-hand side using the following result, which can be
deduced from the proof of Theorem 2.4.3 of \cite{Mi}.
Let $f$ be a modular form of weight $k$ for a Fuchsian group $\Gamma$ commensurate with $SL_2 (\ZZ)$.
Then
$${\rm deg}(f) = \frac{k}{12} \cdot
\frac{[SL_2 ( \ZZ) : \Gamma \cap SL_2 ( \ZZ) ]}{[\Gamma: \Gamma \cap SL_2 ( \ZZ) ]} ~.
$$
This determines ${\rm deg} (\Theta_{\CN} (z))$. ~~~$\bsq$

From now on let $N$ be a fixed number from \eqn{Eq4}.
Define
$$g_1 (z) = \Theta_{\CN} (z) ~,$$
and let $g_2 (z)$ be a modular function for $\frac{1}{2} \Gamma_0 (4N)^+$ with divisor ${\bf \infty} - {\bf 1}$ (which exists since
$\frac{1}{2} \Gamma_0 (4N)^+$ has genus 0).
To be precise, let
$$\eta^{(N)} (z) = \prod_{d|N} \eta (dz ) ~.$$
Then (cf. \cite[Table 3]{Con17}) if $N$ is odd we take
$$g_2 (z) = \left\{
\frac{\eta^{(N)} \left( \frac{z}{2} \right) \eta^{(N)} (2z)}{\eta^{(N)} (z)^2} \right\}^{D_N/ \dim \, \CN} ~,
$$
and if $N$ is even we take
$$g_2 (z) =
\left\{
\frac{\eta^{(N/2)} \left( \frac{z}{2} \right) \eta^{(N/2)} (4z)}{\eta^{(N/2)} (z) \eta^{(N/2)} (2z)} \right\}^{D_N/ \dim \, \CN}~.
$$
\begin{theorem}\label{C2}
Any element $f(z)$ of
$$\sM_{kw^{(N)}} \left( \frac{1}{2} \Gamma_0 (4N)^+ , ( \chi^{(N)} )^k \right)$$
can be written uniquely as
\beql{EqQ6a}
f(z) = g_1 (z)^k \sum_{i=0}^{\lfloor k ~{\rm ord}_1 (g_1) \rfloor} c_i g_2 (z)^i ~.
\eeq
For a cusp form, $c_0 =0$, and if $k ~{\rm ord}_1 (g_1)$ is an integer then that coefficient must also be zero.
\end{theorem}

\paragraph{Proof.}
$f(z)/g_1(z)^k$ is a modular function for ${1\over 2}\Gamma_0(4N)^{+}$
with the trivial character, and therefore can be written as a rational
function in $g_2(z)$.  But, since $f(z)$ has no poles, the only pole
of $f(z)/g_1(z)^k$ is at the cusp class 1, which is also the pole of
$g_2(z)$.  It follows that $f(z)/g_1(z)^k$ is a polynomial in
$g_2(z)$.  The remaining statements follow by considering the order of
$f(z)$ at the two cusp classes.~~~$\bsq$

There is an expression similar to \eqn{EqQ6a} for the theta series of the $\emptyset$-shadow.
Let $s=D_N/ \dim \, \CN$, and set
\begin{eqnarray*}
s_1 (z) & = & \left(\frac{n}{\sqrt{\sqrt{N}z}} \right)^{\dim \CN} g_1 \left( 1- \frac{1}{N z} \right) =
\prod_{d|N} \theta_2 (dz) =
2^{d(N)} \frac{\eta^{(N)} (2z)^2}{\eta^{(N)} (z)} ~, \\
s_2 (z) & = & g_2 \left( 1- \frac{1}{Nz} \right) = - 2^{-D_N/2}
\left\{ \frac{\eta^{(N)} (z)}{\eta^{(N)} (2z)} \right\}^s ~.
\end{eqnarray*}

\begin{coro}\label{C3}
If $\La$ is a strongly $N$-modular lattice that is rationally equivalent to
$(\CN )^k$ then its theta series can be written in the form \eqn{EqQ6a}, and its $\emptyset$-shadow $S$ has theta series
\beql{EqQ6b}
\Theta_S(z) ~=~ s_1 (z)^k \sum_i c_i s_2 (z)^i ~.
\eeq
\end{coro}
\paragraph{Proof.}
This follows from Corollary~\ref{Co3}, Theorem~\ref{C2} and
Equation~\eqn{Eq14a}.~~~$\bsq$

The proofs of Theorems~\ref{th1} and \ref{th2} use only the nonnegativity of the coefficients of certain theta series.
In some cases stronger bounds may be obtained by using the facts that the coefficients must also be integers,
or, more precisely, that $\Theta_\La$ and $\Theta_S$ must have nonnegative
integer coefficients and satisfy $\Theta_\La \equiv 1$ $(\bmod~2)$ and $\Theta_S \equiv 0$ or 1 $(\bmod~2)$;
and if $\La$ is odd with minimal norm $\mu$ then
\begin{eqnarray*}
&& \# \{ v \in S_{\emptyset}( \La ) : ~ v \cdot v  < \mu /4 \} = 0 ~, \\
&& \# \{ v \in S_{\emptyset}( \La ): ~ v \cdot v < \mu / 2 \} \le 2 ~.
\end{eqnarray*}

For example, let us prove that there is no 14-dimensional 3-modular lattice meeting the bound of Theorem \ref{th2} (and satisfying the hypothesis of that theorem).
For such a lattice, Corollary~\ref{C3} would imply that
\begin{eqnarray*}
\Theta_\La  & = & g_1^7 (c_0 + c_1 g_2 + c_2 g_2^2 + c_3 g_2^3 ) \\
& = & 1 + O(q^4) ~, \\
\Theta_S & = & s_1^7 (c_0 + c_1 s_2 + c_2 s_2^2 + c_2 s_2^3 ) ~.
\end{eqnarray*}
From the first equation we find that
$c_0 =1$, $c_1 = -14$, $c_2 = 28$, $c_3 = - 56$, so $\Theta_\La = 1+ 602 q^4 + 1344q^5 + 4032 q^6 + \cdots$, and
then $\Theta_S = \frac{7}{2} q + \frac{147}{2} q^3 + \cdots$, which is impossible.

The nonexistence results for $N=2$ and 3 mentioned at the beginning of Section 2 (and further results given in the table in \cite{SN}) were obtained in this way.

\section{The proofs of Theorems~\ref{th1} and \ref{th2}}
\hsp
We begin by stating a series of identities that relate
$g_1$, $g_2$, $s_1$ and $s_2$.
(We include more than are needed, because of their intrinsic interest.)
For $N$ odd, we have
\beql{ED1}
g_2 (z)^2 g_1 (z)^s = \eta^{(N)} (z)^s ~, 
\eeq
\beql{ED2}
g_2 (z) g_1 (z)^s = - \eta^{(N)} \left( \frac{z+1}{2} \right)^s ~,
\eeq
\beql{ED3}
s_2 (z)^2 s_1 (z)^s = \eta^{(N)} (z)^s ~,
\eeq
\beql{ED4}
s_2 (z) s_1 (z)^s = - 2 ^{D_N /2}~ \eta^{(N)} (2z)^s~,
\eeq
\beql{ED5}
g_2 (z) g_2 (z+1) s_2 (z) = 2^{-D_N /2} ~,
\eeq
\beql{ED6}
\frac{1}{g_2 (z)} + \frac{1}{g_2 (z+1)} + \frac{1}{s_2 (z)} = 2s ~,
\eeq
\beql{ED7}
g_1 (z) g_1 (z+1) s_1 (z) = 2^{\dim \, \CN} \eta^{(N)} (z)^3 ~,
\eeq
\beql{ED8}
g_1 (z)^{s/2} - g_1 (z+1)^{s/2} - s_1 (z)^{s/2} = 2s \eta^{(N)} (z)^{s/2} ~.
\eeq

For $N=2$, $s_1$ and $s_2$ are given by
\begin{eqnarray*}
s_1 (z) & = & \frac{2 \eta (z)^5 \eta (4z)^2}{\eta \left( \frac{z}{2} \right)^2 \eta (2z)^3} ~, \\
s_2 (z) & = & - \frac{1}{16}
\frac{\eta \left( \frac{z}{2} \right)^8 \eta (2z)^{16}}{\eta (z)^{16} \eta (4z)^8} ~.
\end{eqnarray*}
Then we have
\begin{eqnarray}\label{ED9}
&& g_1 (z)^8 g_2 (z)^2 = \eta (z)^8 \eta (2z)^8 ~, \\
\label{ED10}
&& s_1 (z)^4 s_2 (z)^2 = \frac{1}{16} \theta_4 (z)^4 \theta_3 (2z)^4 ~, \\
\label{ED11}
&& s_1 (z)^8 s_2 (z)^2 = \eta (z)^8 \eta (2z)^8 ~, \\
\label{ED12}
&& g_2 (z) s_2 (z+1) = g_2 (z+1) s_2 (z) = \frac{1}{16} ~,
\\
\label{ED13}
&& g_2 (z) g_2 (z+1) s_2 (z) s_2 (z+1) = \frac{1}{256} ~, \\
\label{ED14}
&& \frac{1}{g_2 (z)} + \frac{1}{g_2 (z+1)} + \frac{1}{s_2 (z)} +
\frac{1}{s_2 (z+1)} = 16 ~.
\end{eqnarray}

To show \eqn{ED5}, for example, we observe that $f(z) = g_2 (z) g_2 (z+1) s_2 (z)$ is invariant under $\Gamma_0 (N)^+$, which is transitive on cusps,
and so the order of $f(z)$ at every cusp is the same.
On the other hand every pole and zero of $f(z)$ is at a cusp, and since $f(z)$ is a modular function the number of zeros must equal the number of poles.
Therefore $f(z)$ has no zeros or poles, and must be constant.
We leave the proofs of the other identities to the reader.
\paragraph{Proof of Theorem \ref{th1}.}
Let $\La$ be a unimodular lattice of minimal norm $\mu$ in dimension $n= 8t + o = 24m -l$,
where $0 \le o \le 7$, $1 \le l \le 24$.
We must show that $\mu \le 2m$ except when $m=l=1$.
From Corollary~\ref{C3},
\beql{EQ14a}
\Theta_\La = g_1^n \sum_{i=0}^t c_i g_2^i = \sum_{j=0}^\infty a_j q^j , \quad \mbox{(say)} ~,
\eeq
and the theta series of the $\emptyset$-shadow of $\La$ is
\beql{EQ14b}
\Theta_S = s_1^n \sum_{i=0}^t c_i s_2^i = q^{\frac{o}{4}} \sum_{j=0}^\infty b_j q^{2j} , \quad \mbox{(say)} ~.
\eeq
Suppose, seeking a contradiction, that $\mu \ge 2m+1$.
Then $\Theta_\La = 1+ O (q^{2m+1})$.
This determines $c_i$ for $0 \le i \le 2m$.
In particular, as we show below,
$c_{2m} \le 0$,
with equality only when $n=23$.
On the other hand, we will also write $c_{2m}$ as a linear
combination of $b_j$ for $0 \le j \le t-2m$ with nonnegative coefficients,
and thus $c_{2m} \ge 0$, which is a contradiction unless $n=23$.

To compute $c_{2m}$ we divide both sides of
\eqn{EQ14a} by $g_1^n$ to obtain
$$g_1^{-n} + O(q^{2m+1} ) = \sum_{i=0}^\infty c_i g_2^i ~,$$
where we adopt the convention that
$c_i =0$ for $i > t$.
From the B\"{u}rmann-Lagrange theorem \cite{MS1414} we deduce
$$c_i = - \frac{n}{i}~ \mbox{coefft. of $q^i$ in} ~
q~g'_1~ g_1^{-n-1} \left(
\frac{q}{g_2} \right)^i ~,
$$
for $0 \le i \le 2m$.
For $i=2m$ this simplifies to
$$- \frac{24m-l}{2m} ~\mbox{coefft. of $q^{2m}$ in
$q~ g'_1~ g_1^{l-1}~ q^{2m}~ \eta^{-24m}$} ~,
$$
using \eqn{ED1}.
Now $g'_1 g_1^{l-1}$ (the derivative of a theta series) has nonnegative coefficients,
and $q^{2m} \eta^{-24m}$ has nonnegative coefficients and positive coefficients at even powers of $q$.
So as long as $q g'_1 g_1^{l-1}$ has a nonzero coefficient of even degree
$\le 2m$, it follows that $c_{2m} < 0$.
Since $g_1$ has a linear term and $g'_1$ has a cubic term, the only way $c_{2m}$ can equal zero is if $m=l=1$, i.e. if $n=23$.

On the other hand, from \eqn{EQ14b} we have
$$
\sum_{i=0}^t c_{t-i} s_2^{-i} = s_1^{-n} s_2^{-t} q^{o/4} \sum_{j=0}^\infty b_j q^{2j}$$
and thus
$$c_i = \sum_{j=0}^{t-j} \beta_{i,j} b_j$$
where
$$s_1^{-n} s_2^{-t} q^{o/4+2j} = \sum_{i=0}^\infty \beta_{t-i,j}
s_2^{-i} ~.
$$
Again using B\"{u}rmann-Lagrange we find
\beql{EQ19a}
\beta_{i,j} = - ~\mbox{coefft. of $q^{2t-2i-2j}$ in} ~
2~\frac{qs'_2}{s_2}~q^{2t-2i + o/4}~s_2^{-i}~s_1^{-n} ~.
\eeq
From the product expansion for $s_2$ we immediately deduce that all coefficients
of $q s'_2 / s_2$ are nonpositive.
From \eqn{ED4}
(with $s=24$) and the fact that $s_1$ has nonnegative coefficients, the
remaining factor in \eqn{EQ19a} has nonnegative coefficients as long as $24i \ge n$.
In particular, this is certainly true for $i=2m$, and thus $\beta_{2m,j} \ge 0$,
which produces the desired contradiction.~~~$\bsq$
\paragraph{Proof of Theorem \ref{th2}.}
Theorem \ref{th1} covers the case $N=1$,
and the other cases when $N$ is odd are analogous.
The proof for $N=2$ is given below, the remaining cases 6 and 14 again being analogous.
We begin with a lemma.
\begin{lemma}\label{LN2}
If $N=2$ then (i)~all coefficients of $s_1^{-2a} s_2^{-2b}$ are nonnegative
whenever $2b \le a \le 4b$,
and (ii)~all coefficients of the
logarithmic derivative of $q^c s_1^{-2a} s_2^{-2b}$ are nonnegative whenever
$2b \le a \le \min \{2b+c, 4b\}$.
\end{lemma}
\paragraph{Proof.}
We write
$$q^c s_1^{-2a} s_2^{-2b} = q^{2b+ c-a} (q^2 s_1^{-8} s_2^{-2} )^{a/2-b} (s_1^{-4} s_2^{-2} )^{2b-a/2} ~,
$$
in which
the exponents $2b-a/2$, $a/2-b$ and $2b+c-a$ are positive by hypothesis, and consider each factor separately.
First, $q$ and its logarithmic derivative $(q^{-1})$ are both nonnegative.
The other two terms may be expanded as
\beql{EQ21a}
\frac{q^2}{s_1^8 s_2^2} = \frac{q^2}{(\eta^{(2)} )^8} =
\prod_{m=1}^\infty (1-q^m)^{-8} (1-q^{2m})^{-8}
\eeq
and, using \eqn{ED10},
\beql{EQ21b}
s_1^{-4} s_2^{-2} = 16 \prod_{m=1}^\infty
\{ (1+q^{2m-1} ) (1-q^{2m-1} )^{-2} \}^4
(1-q^{4m} )^{-8} (1-q^{8m-4} )^{-8} ~.
\eeq
Each factor of \eqn{EQ21a} and \eqn{EQ21b} has nonnegative coefficients and
nonnegative logarithmic derivative.~~~$\bsq$
\paragraph{Proof of Theorem \ref{th2} for $N=2$.}
Let $\La$ be a 2-modular lattice of minimal norm $\mu$ in dimension $n= 4t + o = 16m -2l$,
where $o=0$ or 2, $1 \le l \le 8$.
Then
\begin{eqnarray*}
\Theta_\La & = & g_1^{n/2} \sum_{i=0}^t c_i g_2^i = \sum_{j=0}^\infty a_j q^j ~, \\
\Theta_S & = & s_1^{n/2} \sum_{i=0}^t c_i s_2^i = q^{o/2}
\sum_{j=0}^\infty b_j q^j ~,
\end{eqnarray*}
say.
That $c_{2m} \le 0$ follows as in the proof of Theorem~\ref{th1}, but
since $g'_1$ now has a linear term there is no exception and $c_{2m} < 0$.

On the other hand, defining $\beta_{i,j}$ as before, we find
$$
\beta_{2m,j} = \frac{1}{2m} ~\mbox{coefft. $q^{t-2m-j}$ in} 
~ \sL \{ q^{j+ o/2} s_1^{-n/2} s_2^{-t} \} q^{n/4 -2m} s_1^{-n/2} s_2^{-2m} ~,
$$
where $\sL$ denotes the logarithmic derivative.
By the lemma, the first factor has nonnegative coefficients when $0 \le o/2 \le t$, and the second factor has nonnegative coefficients when $8m \le n \le 16m$.
This proves the desired result for $n \ge 8$.
The remaining three cases, $n=2,4$ and 6, can be checked directly.~~~$\bsq$

We end this section with 
an analogue of Theorem~\ref{th1} for codes over $\ZZ/4 \ZZ$.
This is a generalization of a bound established by
Bonnecaze et~al. \cite{BSBM97} for self-dual codes over
$\ZZ/4 \ZZ$ in which all Euclidean norms are divisible by 8.
\begin{theorem}\label{thZ41}
Suppose $C$ is a self-dual code over $\ZZ /4 \ZZ$ of length $n$.
The minimal Euclidean norm of $C$ is at most
$$8 \left[ \frac{n}{24} \right] + 8 ~,
$$
unless $n \equiv 23$ $(\bmod ~24)$ when the bound must be increased by $4$.
\end{theorem}
\paragraph{Proof.}
As in \cite{BSBM97} we construct a unimodular lattice $\La$ from $C$ using
``Construction A''.
$\La$ has theta series
$$\theta_3 (4z)^n + O(q^{\mu /4} ) ~,
$$
where $\mu$ is the minimal Euclidean norm of $C$.
The argument used to prove Theorem~\ref{th1} now establishes the desired result.
The identity
$$
q ~ ( \theta_3(4z)~\theta_3'(z) ~-~ \theta_3'(4z)~\theta_3(z) ) ~=~ \eta(2z)^6
$$
is needed.  ~~~$\bsq$

\section{Genera not covered by Theorem \ref{th2}}
\hsp
We can say less about the minimal norm $\mu$ of a strongly $N$-modular lattice $\La$ not rationally equivalent to a direct sum of copies of $\CN$.
If $N=1$ or 2, local considerations show that no such lattices exist,
nor can they exist for $N=7$ or 23,
although then a more involved argument appears to be needed (see the Appendix).

For $N=3,5,6$ and 11, the numerical evidence suggests that
\beql{Eq5a}
\mu \le 2 \left[ \frac{n- \dim \, \CN}{D_N} \right] +2 ~.
\eeq
For $N=3$ and 11 we further conjecture that if equality holds and $n \equiv 2$ $(\bmod~D_N$) then $\La$ must be even.
These conjectures have been verified for $n \le 56$ for $N=3,5,6$ and for $n \le 32$ for $N=11$.

For $N \le 14$ we conjecture that
$$\mu \le 2 \left[ \frac{n}{4} \right ] +2 , \quad n \neq 2 ~;
$$
this has been verified for $n \le 30$.
For $N=15$ there is no obvious pattern.
In the critical dimension 4, for example, the lattice
defined below in \eqn{EqA5} has minimal norm 4, which actually coincides with the bound of Theorem~\ref{th2}.
\paragraph{Acknowledgement.}
The computer language Magma \cite{Mag1}, \cite{Mag2}, \cite{Mag3}
has been helpful in studying particular lattices, testing for modularity, etc.
\clearpage
\section*{Appendix}
\hsp
For $N=7$ and $23$, every $N$-modular lattice
must satisfy the hypotheses of Theorem \ref{th2}.
This a consequence of the
following result.

\begin{theorem}\label{thA1}
If $N$ is a positive integer congruent to $7$ mod $8$, then any
even $\{N\}$-modular lattice of even-level $4N$ has oddity $0$.
\end{theorem}

We begin with two lemmas.  When we say that a power series reduces to 1 mod
2, this includes the assertion that the coefficients are algebraic
integers, and the corresponding number field is unramified at 2.

\begin{lemma}\label{LemA2}
Let $\Lambda$ be an even lattice of even-level $2^k N$, with $N$ odd.  Then
for any element $t$ of $\Gamma_0(2^k)$, there exists a constant $C$ such
that $C \Theta_\Lambda |_t$ reduces to $1$ mod $2$.
\end{lemma}

\paragraph{Proof.}
The analysis of \cite{Kit80} can be extended to show that there exists
a function $T(v)$ on $\Lambda^*$
such that
$$
\Theta_\Lambda |_t = \sum_{v\in\Lambda^*} T(v) q^{v \cdot v}.
$$
Moreover $T(v)=T(-v)$ and $T(v)/T(0)$ is either 0 or a root of unity of odd
order.  Taking $C=T(0)^{-1}$, the result follows immediately.~~~$\bsq$

\begin{lemma}\label{LemA3}
Let $g$ be a modular function for $\Gamma(2)$.  If
all poles of $g$ occur at the cusp $1$, and the expansion of $g$ at
$\infty$ reduces to $1$ mod $2$, then all zeros of $g$ occur at points $z$ such
that $16\lambda(z)^{-1}$ is an algebraic integer with even norm, where
$\lambda(z)=(\theta_2(z)/\theta_3(z))^4$.
\end{lemma}

\paragraph{Proof.}
Let $\ell(z)$ be the function $\lambda(z)/16$.  Then
$\ell(z)$ has integer coefficients, with leading coefficient equal to 1.
Since the unique pole of $\ell$ is at $1$, $g$ can be expressed as
a polynomial $p$ in $\ell$, with algebraic integer coefficients.
Clearly, then, $g$ reduces to 1 mod 2 just when $p$ reduces to 1 mod
2.  But this implies that all roots of $p$ have 2-adic valuation greater
than 1, which is the desired result.~~~$\bsq$

\paragraph{Proof of Theorem \ref{thA1}.}
Let $\Lambda$ be an $n$-dimensional even
$\{N\}$-modular lattice of even-level $4N$, with theta series
$\Theta_\Lambda(z)$, and let $K$ be the $\{4,N\}$-modular lattice
$
\sqrt{2} {\Bbb Z}\times \sqrt{2N} {\Bbb Z},
$
with theta series
$\Theta_K(z)=\theta_3(2z)\theta_3(2Nz)$.
Then
$
f(z)=\Theta_\Lambda(z)/\Theta_K(z)^{n/2}
$
is a modular function for $\Gamma_0(4N)$ (with trivial character,
since $\dim \La$ is even and thus $\det_2 \La$ is a square).
Furthermore, since $\theta_3(2z)\theta_3(2Nz)$ has zeros only at cusps, it
follows that $f$ has poles only at cusps.

If $\Lambda$ had oddity 4 (the only other possibility), then $f$
would satisfy the relation
$$
f|_{W_N} = -f,
$$
since both $\Lambda$ and $K^{n/2}$ are $\{N\}$-modular, and $K$
has oddity 0.  As a consequence, $f$ has at least one zero at every
point of $\Gamma_0(4N)$ fixed by $W_N$.  Also, since $f$ is the ratio
of two theta series, its expansion around $\infty$ has integer coefficients,
and reduces to 1 mod 2.

Let $T$ be a set of (right) representatives for $\Gamma_0(4)/\Gamma_0(4N)$.
Then
$$
g=\prod_{t\in T} f|_t
$$
is a modular function for $\Gamma_0(4)$ ($g$ can also be thought of as the
norm of $f$ from $\Gamma_0(4N)$ to $\Gamma_0(4)$).  Moreover, up to a
constant factor, the expansion of $g$ around $\infty$ reduces to 1 mod 2,
since by Lemma \ref{LemA2} the same is true for each $f|_t$.

Since $g(z)$ is invariant under $\Gamma_0(4)$, $g(z/2)$ is invariant under
${1\over 2}\Gamma_0(4)=\Gamma(2)$.  Moreover, since the poles of $g(z)$ are
at the cusps, and neither $\infty$ nor $0$ are poles of $g(z)$, it follows
that the only pole (which may, of course, be a multiple pole) of $g(z)$ is
at the cusp $1$.

To finish the proof, we invoke Lemma \ref{LemA3}.  We obtain a contradiction if
we can demonstrate the existence of some point $z$ of $\Gamma_0(4N)$ fixed
by $W_N$ such that $16\lambda(z/2)$ has odd norm.  Let $x$ be the point
$(1+\sqrt{-N})/4$ of $\Gamma_0(N)$.  If
$x'$ is any image of $x$ in $\Gamma_0(N)\cap \Gamma(2)$, then $2x'$ is a point
of $\Gamma_0(4N)$ fixed by $W_N$.  Now $j(x)$
has odd norm. (The elliptic curve corresponding to $x$ has complex
multiplication by an order of ${\Bbb Q}(\sqrt{-N})$ of odd conductor.
A CM curve always has integral $j$-invariant, so has good
reduction over a suitable extension of ${\Bbb Q}$.  The reduction mod 2
cannot be supersingular, so the $j$-invariant cannot reduce to 0.) We have
$$
j={(\ell^2-16\ell+256)^3\over \ell^2(\ell-16)^2}\equiv \ell^2\pmod 2,
$$
so for two of the six images, $\ell$ must have odd norm.~~~$\bsq$

\paragraph{Remarks.}
(1) The assumption that $N$ is congruent to 7 mod 8 is
critical; for $N$ congruent to 1 mod 4, there are no points fixed by $W_N$
other than cusps, while for $N$ congruent to 3 mod 8, the points fixed by
$W_N$ correspond to curves with supersingular reduction mod 2.  (2) The
hypothesis that the even-level be $4N$ can be relaxed to say that the
even-level is $4MN$, where $M$ is an odd integer,
relatively prime to $N$, such that $-N$ has a square root $\bmod~M$
and $\det_{\Pi (M)} \La$ is a square;
in that case, the conclusion is that
$\gamma_{\Pi(4M)}(\Lambda)=1$.  The existence of a square root of $-N$ is
necessary to allow the existence of suitable CM curves.

\begin{coro}\label{CoA4}
If $N$ is an integer congruent to $7$ mod $8$, then any $N$-modular 
lattice has oddity $0$, as does any
$\{2,N\}$-modular lattice of dimension a multiple of $4$.
\end{coro}

\paragraph{Proof.}
If $\Lambda$ is an $N$-modular lattice, then $\sqrt{2}\Lambda$
is $\{4,N\}$-modular, and has the same oddity (since $\det{\Lambda}$ is
1 or 7 mod 8).  Therefore $\sqrt{2}\Lambda$ satisfies the hypotheses
of Theorem \ref{thA1}, and must have oddity 0.

Similarly, a $\{2,N\}$-modular even lattice has oddity 0.  If
$\Lambda$ is a $\{2,N\}$-modular odd lattice, then the even neighbor
of $\Lambda$ (recall Theorem \ref{ThEN}) is a $\{2,N\}$-modular even lattice
with the same oddity.~~~$\bsq$

The following is immediate:

\begin{coro}\label{CoA5}
All $p$-modular lattices, for $p$ prime and
congruent to $7$ mod $8$, are rationally equivalent to the direct sum of some
number of copies of $C^{(p)}$.  A strongly $14$-modular lattice must be
rationally equivalent to the direct sum of some number of copies of
$\left( \begin{array}{cc} 3&1\\1&5 \end{array}\right)$.
A strongly $15$-modular lattice is rationally equivalent to
the direct sum of some number of copies of $C^{(15)}$, possibly together
with a copy of
\beql{EqA5}
\left(
\begin{array}{cccc}
4&0&2&1\\
0&4&1&2\\
2&1&5&1\\
1&2&1&5
\end{array}
\right)~.
\eeq
\end{coro}

\end{document}